\documentstyle[12pt]{article}
\begin{document}

\newcommand{\cF}{{\cal F}}
\newcommand{\cL}{{\cal L}}
\newcommand{\cC}{{\cal C}}
\newcommand{\cQ}{{\cal Q}}
\newcommand{\cD}{{\cal D}}
\newcommand{\eps}{\epsilon}
\newcommand{\exist}{\exists}
\newcommand{\x}{{\bf x}}
\newcommand{\y}{{\bf y}}
\newcommand{\z}{{\bf z}}
\newcommand{\X}{{\bf x}}
\newcommand{\Y}{{\bf y}}
\newcommand{\D}{{\bf d}}
\newcommand{\glceil}{\left\lceil}
\newcommand{\grceil}{\right\rceil}
\newcommand{\glfloor}{\left\lfloor}
\newcommand{\grfloor}{\right\rfloor}
\newcommand{\qed}{\hfill$\Box$}
\newcommand{\proof}{\smallskip\par\noindent {\bf Proof: }}
\newcommand{\aleq}{\stackrel{\sim}{<}}
\newcommand{\ef}{{\stackrel{\rm \Delta}{=}}}

\newtheorem{teor}{Theorem}
\newtheorem{defi}{Definition}
\newtheorem{lem}{Lemma}
\newtheorem{claim}{Claim}
\newtheorem{cor}{Corollary}
\newtheorem{pro}{Problem}
\newtheorem{prop}{Property}

\baselineskip=0.7truecm
\begin{titlepage}
\title{A splitting lemma}
\author{{\bf G. Greco} \\
Objectway spa\\
\\
e-mail:{\tt giampaolo.greco@objectway.it }}
\thanks{Via F.Domiziano 10, 00147 Roma.}

\maketitle

\begin{abstract}

In this paper, we study the relations between the numerical structure of the optimal solutions
of a convex programming problem defined on the edge set of a simple graph and the
{\it stability number} (i.e. the maximum size of a subset of pairwise non-adjacent vertices) of
the graph.  Our analysis shows that the stability number of every graph $G$ can be decomposed in
the sum of the stability number of a subgraph containing a perfect $2$-matching (i.e. a system of
vertex-disjoint odd-cycles and edges covering the vertex-set) plus a term computable in polynomial
time. As a consequence, it is possible to bound from above and below the stability number in terms of the
{\it matching number} of a subgraph having a perfect $2$-matching and other quantities computable
in polynomial time.  Our results are closely related to those by Lorentzen~\cite{Lor}, Balinsky and
Spielberg~\cite{BS}, and Pulleyblank~\cite{Pul} on the linear relaxation of the
{\it vertex-cover problem}. Moreover, The convex programming problem involved has important
applications in information theory and extremal set theory where, as a graph capacity formula, has
been used to answer some longstanding open questions~(see \cite{CKS} and~\cite{GKV1}).

\bigskip

{\flushleft {\bf keywords. }}matching, $2$-matching, stability number,
packing, covering, entropy, graph capacity.

\end{abstract}

\end{titlepage}

\section{Terminology and notation}

Given any two positive reals $0<p,q<1$ we define the function
$$
\hbar(p,q) \ef (p+q)h \left ( \frac{p}{p+q} \right ),
$$
where
$$
h(x)=-x\log x-(1-x)\log (1-x), \qquad (0<x<1)
$$
is the {\it binary entropy} and (here and in the sequel) log's are to the base $2$.

A $stable$ set in a simple graph $G$ is a set of vertices that does not
contain any edge. The size of a maximum stable set in $G$ is the
{\it stability number} of $G$ and it is denoted by $\alpha(G)$. A set of vertices
is a {\it vertex cover} of $G$ if each edge has at least one endpoint in
the set. The minimum size of a vertex cover is the {\it cover number} of
$G$ and it is denote by $\tau(G)$. These two invariants are closely related
by the Gallai identity:
$$
\alpha(G)+\tau(G)=n \qquad (n=|V(G)|).
$$
The {\it matching number} of $G$, denoted by $\nu(G)$, is the maximum number of
disjoint edges contained in the graph.

If $X$ is a set of vertices of $G$, we denote by $\bar{X}=V(G) \setminus
X$ the complementary of $X$. Moreover $\Gamma X$ is the set of
vertices of $\bar{X}$ adjacent to some vertex in $X$. Finally, if
$X \subseteq V(G)$ and $F=(V(F),E(F))$ is a subgraph of $G$:
$$
\bar{X}(F)=\bar{X} \cap V(F) \:\:\hbox{ and }\:\: X(F)=X \cap V(F).
$$

\section{Introduction}
\message{Introduction}

In this paper we study the numerical structure of the optimal solutions of
the following convex programming problem defined on the edge set of a graph. Let $G=(V(G),E(G))$ be
a simple graph, $P$ a distribution of probability defined on $V(G)$ and set:
$$
l(G,P)=\min_{\{x,y\} \in E(G)}\hbar(P(x),P(y)).
$$
We define the $conjunctive$ $capacity$ of $G$ as:
\begin{equation}
\label{cap}
\Theta(G)=\max_{P}l(G,P).
\end{equation}
Note that, being $\hbar$ concave, problem~(\ref{cap}) can be reduced to a convex
programming problem. A distribution $P$ is $G$-balanced if
it achieves the maximum in~(\ref{cap}). Let us define the $t$-th
power of $G$ as the graph $G^t=(V(G)^t,E(G^t))$ such that
$\{(x_1,\ldots,x_t),(y_1,\ldots,y_t)\} \in E(G^t)$ if for every edge $e \in E(G)$
there exists a position $1\leq i \leq t$ such that $\{x_i,y_i\}=e$. In~\cite{CKS}
and~\cite{GKV1} the authors show that $\Theta(G)$ is the asymptotic exponent of the
clique number (i.e. the size of the largest complete subgraph of $G$) of the powers
of $G$:
$$
\Theta(G)=\lim_{t \rightarrow +\infty}\frac{1}{t}\log\omega(G^t).
$$
This result has been used to answer a long-standing open question on the asymptotics
of the maximum number of qualitatively independent partitions in the sense of
R\'enyi~\cite{Re}. We point out that in these papers the conjunctive capacity of graphs
is considered as a particular case of the Sperner capacity of a family of directed graphs.
For the applications in information theory see for example~\cite{CKS}.

By considering the uniform distribution on the vertex set of the graph, one easily see
that for every graph with $n$ vertices, $\Theta(G) \geq 2/n$.
A $2$-matching is a vector ${\bf x}=(x_e:~e \in E(G))$ with components
$0$, $1$ or $1/2$, such that for every $x \in V(G)$ the sum of the weights to the
edges incident in $x$ is at most $1$. A $2$-matching is maximum if the overall sum of the
weights assigned to the edges is maximum. A $2$-matching is $perfect$ if every vertex has
some incident edge with nonzero weight. It is easy to see that a graph has a perfect $2$-matching
if and only if it contains a system of vertex-disjoint odd cycles and edges covering the
vertex-set (for more on $2$-matchings and related problems see~\cite{LP}). In~\cite{Gre95} we
show the following
\begin{teor}
\label{G1}
For every simple graph $G$ without isolated vertices,
$$
\Theta(G)=\frac{2}{n}, \qquad (n=|V(G)|)
$$
if and only if $G$ has a perfect $2$-matching.
\end{teor}

Note that an easy corollary is that if the uniform distribution is
$G$-balanced then this it is also the unique optimal solution to~(\ref{cap}).
We recall also the following characterization for graph having a perfect
$2$-matching~\cite{Tut}:

{ \flushleft {\bf Theorem T }
$G$ has a perfect $2$-matching if and only if for every stable set $X \subseteq V(G)$:
$$
\frac{|\Gamma X|}{|X|} \geq 1.
$$}

In the sequel $G$ is a graph without isolated vertices. A vertex is
$critical$ if its deletion strictly decreases the stability number of the
graph. It is easy to see that a vertex is critical if and only if it
belongs to every stable set of maximum size. If $P$ is a probability
distribution on $V(G)$, the vertex $x$ is $P-critical$ if, for some
$y \in \Gamma x$, $P(x) < P(y)$. In the next section we prove the following
\begin{lem}{\bf (Splitting Lemma)}
\label{main}
For any graph $G$ and $G$-balanced distribution $P$,
\begin{enumerate}
\item All the $P$-critical vertices are critical,
\item If $X$ is the set of the $P$-critical vertices, then the subgraph
of $G$:
$$
F=G-(X \cup \Gamma X)
$$
has a perfect $2$-matching.
\end{enumerate}
\end{lem}

Note that an immediate consequence of the Splitting Lemma is the
following, already known, result (for example, see~\cite{Ber}):

\begin{cor}
If $G$ has no critical points then it has a perfect $2$-matching.
\end{cor}
\proof
By the hypothesis and the Splitting Lemma it follows that every $G$-balanced distribution
has an empty set of $P$-critical vertices and $G$ has a perfect $2$-matching.
\qed

\medskip

The number of the $P$-critical vertices is computable in polynomial time. So, it is interesting to
investigate its relations with the stability number of $G$:

\begin{teor}
\label{alpha}
Let $G$ be a graph and $X$ the set of $P$-critical vertices for a
$G$-balanced distribution $P$. Then
$$
\alpha(G)=|X|+\alpha(F),
$$
where $F=G-(X \cup \Gamma X)$ has a perfect $2$-matching.
\end{teor}

So, the stability number of every graph can be expressed as the sum of
the stability number of a graph with a perfect $2$-matching plus some
quantity computable in polynomial time. Now, the stability number of a graph
with a perfect $2$-matching can be bounded from above and below in terms
of the matching number of the graph. Indeed, by observing that the set of
vertices non covered by a maximal matching of $G$ is a stable set, one
gets the general lower-bound:
$$
\alpha(G) \geq n-2\nu(G).
$$
On the other hand, if $G$ has a perfect $2$-matching and $X$
is any maximum stable set in $G$ then, by Theorem~T
$$
|\Gamma Y| \geq |Y|, \:\:\:\hbox{for every $Y \subseteq X$}.
$$
By Hall's Theorem, $G$ contains a matching covering every vertex in $X$, and
$$
\alpha(G) \leq \nu(G).
$$
It follows
\begin{cor}
\label{co2}
Let $G$ be a graph, $P$ a $G$-balanced distribution, $X$ the set of
$P$-critical points, and $F=G-[X \cup \Gamma X]$, then:
$$
|X|+|V(F)|-2\nu(F) \leq \alpha(G) \leq |X|+\nu(F),
$$
In particular
$$
\nu(F)=\frac{|V(F)|}{3} \:\:\Rightarrow\:\: \alpha(G)=|X|+\nu(F).
$$
\end{cor}

{\flushleft {\bf Remark:}} Note that if a graph $F$ has a perfect $2$-matching then
$$
\nu(F) \geq \frac{|V(F)|}{3}.
$$

\bigskip

The set $P$-critical points, $P$ balanced, plays a similar role of the set of vertices with weight
zero in a minimum $2-cover$ of $G$. A {\it fractional vertex cover} is any feasible
solution ${\bf y}=(y_u:\,u \in V(G))$, of the following dual of a linear programming
problem
\begin{equation}
\label{linear-cover}
\left\{\begin{array}{l}
min \:\:{\bf 1}{\bf y} \\
y_u+y_v \geq 1 \:\:\:\:\forall\;\{u,v\} \in E(G) \\
{\bf y} \geq 0
\end{array}\right.
\end{equation}
An optimal solution is a $minimum$ fractional cover. A {\it $2$-cover} of $G$ is a
fractional cover whose components are $0$, $1$ or $1/2$. A $2$-cover is $basic$
if the graph induced in $G$ by the set of vertices with weight $1$ is $not$ bipartite.
Lorentzen~\cite{Lor} and independently Balinsky and Spielberg~\cite{BS} proved that the set of
vertices of the feasible region of problem~(\ref{linear-cover}) coincides with the set of
the basic $2$-covers of the graph. It is possible to prove that the uniform $2$-cover
(i.e. the assignment of weight constantly equal to $1/2$) is $an$ optimal solution to
the minimum fractional cover problem if and only if $G$ has a perfect $2$-matching.
Nevertheless, this does not mean that the uniform fractional cover is the $unique$ optimal
solution. For example consider a complete bipartite graph with color classes of same size.
Having this graph a perfect matching, the uniform fractional cover is optimal.
But another optimal solution is the one having value $0$ on a color class and value $1$ on the
complementary class. This simple example shows that the analogous of the Splitting lemma does not
hold for the set of vertices having weight $0$ in an optimal fractional cover of $G$ (Pulleyblank
in~\cite{Pul} prove that the uniform fractional cover is the $unique$ optimal solution if and only if
for every vertex $v \in V(G)$ the graph $G-\{v\}$ has a perfect $2$-matching).

In the next section we give a proof of the Splitting Lemma.

\section{Proof of the Splitting Lemma}
\message{Proof of the Splitting Lemma}


In~\cite{Gre95} we proved the following three lemmas. In all the
statement $G$ has no isolated vertices. In the first lemma, a {\it line cover}
of $G$ is a set of lines collectively incident with each point of $G$:

\begin{lem}
\label{support}~{\rm \cite{Gre95}}
Let $G=(V(G),E(G))$ be a simple graph and $P$ a $G$-balanced distribution,
then:
$$
{\cal L}(P)=\{ \{x,y\} \in E(G) : \hbar(P(x),P(y))=\Theta(G)  \}
$$
is a line cover of $G$.
\end{lem}

Now, set
$$
e(P)=\{x \in V(G):\;P(x)=P(y) \:\:\hbox{for any $y \in \Gamma X$}\},
$$
and let us denote by $m(P)$ the set of $P$-critical vertices in
$G$

\begin{lem}
\label{stables}~{\rm \cite{Gre95}}
Let $G=(V(G),E(G))$ be a simple and $P$ a $G$-balanced distribution. Then $m(P)$ is a stable
set in $G$ and for every maximal stable $S \supseteq m(P)$, $S \setminus
m(P)$ is a maximal stable in the subgraph induced in $G$ by $e(P)$.
\end{lem}

Let $S$ be a maximal stable set of $G$. Every distribution $P$ such that $m(P)
\subseteq S$ is called $centered$ on $S$. The family of all the distributions centered on $S$
will be denoted by $Cr(S)$. Note that the uniform distribution is centered on every
maximal stable set of $G$.

\begin{lem}
\label{two}~{\rm \cite{Gre95}}
Let $G=(V(G),E(G))$ be a simple graph without isolated vertices and $P$
a $G$-balanced distribution centered on $S$. Then for
every  connected component $F=(V(F),E(F))$ of the graph
$(V(G),{\cal L}(P))$ there exist two reals $q_F \leq p_F$ such that
\begin{equation}
P(v)=\cases{p_F & if $v\in V(F) \cap \bar{S}$ \cr
                    & \cr
            q_F & if $v\in V(F) \cap S.$  \cr}
\end{equation}
\end{lem}

Now, it is interesting to consider our maxmin problem for probability distributions that
assume at most two different values on the vertex set of a
graph. In particular if $q$ and $p$ are these two values with $q \leq
p$, by Lemmas~\ref{stables} and \ref{two} there must exist a maximal stable set $S$ in $G$
such that $P(v)=q$ if $v \in S$ and $P(v)=p$ otherwise. In particular, for those graphs $G$
for which there exists a two valued balanced distribution $P$ we obtain the exact solution
of (\ref{main}). Let $G$ be a graph and $S$ a maximal stable set of $G$ with $|S|=\alpha$.
We write $|\bar{S}|=\tau$. Then the maxmin problem for a two-valued distribution can be
defined as :
\begin{equation}
\label{phi-defi}
\phi(w,\alpha,\tau)=\max_{(q,p) \in D_{w,\alpha,\tau}} \hbar(\:p,q\:)
\end{equation}
\begin{equation}
D_{w,\alpha,\tau}=\{(q,p) \in (0,1]^2 \mid q \leq p \:\hbox{and}\:
q\alpha+p\tau=w \},
\end{equation}
where $w=1$ and $\alpha$, $\tau$ are positive constants. In the following
proofs we will consider the general setting where $0 < w \leq 1$.

It will be convenient to rewrite the above, setting
\begin{equation}
t=\frac{p}{q}=\frac{w-q\alpha}{q\tau}
\end{equation}
and define
\begin{equation}
\label{basic}
\phi(w,\alpha,\tau)=\max_{t \,\geq\, 1} \, z(t,w,\alpha,\tau).
\end{equation}
where
\begin{equation}
z(t,w,\alpha,\tau)=\hbar
\left( \:\frac{wt}{t\tau+\alpha}\,,\,\frac{w}{t\tau+\alpha}\: \right).
\end{equation}

Now, we formulate the two main properties of $\hbar$ that will be used in the
sequel (proof in Appendix~A).

\medskip

\begin{prop}
\label{ip1}
$\hbar(.,.) \in {\cal C}^{(1)}(\,(0,1]^2\,)$ is a symmetric and
strictly increasing function of its arguments.
\end{prop}

\begin{prop}
\label{ip2}
For fixed $w,\alpha,\tau$ the function $z(.,w, \alpha , \tau )$
has a unique absolute point of maximum $t(w,\alpha,\tau) \in [1,+\infty)$.
If $t(w,\alpha,\tau)>1$ then it is also the unique stationary point of
$z(.,w,\alpha,\tau)$ and if $t(w,\alpha,\tau)=1$ then $z$ is a strictly
decreasing function for $t>1$.
\end{prop}


\begin{claim}
\label{disjoint}
For any balanced distribution $P$, we have
\begin{equation}
\label{ep}
e(P)=V(G) \setminus [m(P) \cup \Gamma m(P)].
\end{equation}
\end{claim}
\proof
By Property~\ref{ip1}, it is clear that $e(P) \cap m(P)=\emptyset$. Suppose that for an
$x
\in m(P)$, $e(P) \cap \Gamma x \neq \emptyset$ and fix $y  \in e(P) \cap \Gamma x$.
Now,
let $F$ and $F'$ be the connected components in $(V(G),{\cal L}(P))$ containing $x$
and
$y$ respectively.
By $x \in m(P)$ and $y \in e(P)$ $F \neq F'$ and $\{x,y\} \not\in {\cal L}(P)$.
Hence:
$$
p(F') \geq P(x)=q(F') > q(F)=p(F)=P(y),
$$
and by using Property~\ref{ip1} and Lemma~\ref{two} one gets a contradiction with
$$
\hbar(p(F),q(F))=\hbar(p(F'),q(F')).
$$
Therefore
$$
e(P) \subseteq V(G)\setminus[m(P)\cup\Gamma m(P)].
$$
For the converse, suppose $x \not\in e(P) \cup m(P)$. Then for any $\{x,y\} \in
{\cL}(P)$,
$P(x) > P(y)$ and $x \in \Gamma m(P)$.
\qed

Now, we prove item $2$ in the Splitting Lemma. We will use the
following
\begin{equation}
\label{down}
t(w,\alpha,\tau)=1 \hbox{   iff   } \alpha \leq \tau.
\end{equation}
(see Appendix~A)

{\flushleft {\bf Proof of $2$ in Lemma~\ref{main}:}}
By~(\ref{ep}) it suffices to show that, for every $G$-balanced distribution $P$, the
subgraph induced in $G$ by $e(P)$ has a perfect $2$-matching. We have
$$
e(P)=\bigcup_{F:\,q(F)=p(F)}V(F)
$$
where the union ranges into the family of the components $F$ of $(V(G),{\cL}(P))$
such
that $q(F)=p(F)$. We show that if $q(F)=p(F)$ then $F$ has a perfect $2$-matching. By
Tutte's Theorem we must prove that for every stable set $Y$ in $F$,
$$
|\Gamma_FY| \geq |Y| \qquad (\Gamma_FY=V(F) \cap \Gamma Y).
$$
Suppose the contrary and let us fix
$$
t=t(w,|Y|,|\Gamma_FY|), \qquad (|\Gamma_F Y|<|Y|)
$$
where
$$
w=P(Y \cup \Gamma_FY)=(|Y|+|\Gamma_FY|)q \:\:\:\hbox{and $q=q(F)=p(F)$.}
$$
Note that by $|\Gamma_FY| < |Y|$ and~(\ref{down}) $t>1$.

We replace $P$ with a new probability distribution $P'$ where ${\cL}(P')$ is not a line
cover but
$$
l(G,P') \geq l(G,P).
$$
By Lemma~\ref{support} it follows that $P$ cannot be $G$-balanced. Fix
$$
R=|\Gamma_FY||Y|^{-1}
$$
and
$$
\epsilon = \min_{} \left\{ \frac{q(t-1)}{Rt+1}
\,,\,R^{-1}\min_{F' \in {\cC}}{[q-q(F')]} \,,\,\min_{F' \in
{\cC}}{[p(F')-q]} \right\},
$$
$$
\nu = \epsilon R,
$$
where ${\cC}$ is the family of the components $F'$ of $(V(G),{\cL}(P))$ such that $q(F')
\neq q$. Note that $\eps>0$ (in particular, by Property~\ref{ip1}, $p(F')=q$ implies
$q(F')=q$).

Define $P'$ as:
\begin{equation}
P'(v)=\cases{q+\eps & if $v\in \Gamma_FY$ \cr
                    & \cr
            q-\nu & if $v \in Y$  \cr}
\end{equation}
and $P'(v)=P(v)$ if $v \not\in Y \cup \Gamma_FY$. $P'$ is a probability distribution.
Indeed, from $\eps \leq q(t-1)(Rt+1)^{-1}$ it follows
$$
q-\nu=q-\eps R \geq \frac{q(R+1)}{Rt+1}>0,
$$
and $\eps$ and $\nu$ are fixed so as to leave the total amount of probability of $Y \cup
\Gamma_FY$ unchanged. We prove that the global minimum does not decrease, that is
$l(G,P') \geq l(G,P)$.

\medskip

{\bf Case 1: }Edges $\{ x,y \} $ such that one endpoint $x$ belongs to $\Gamma_F Y$. If
$y
\not\in Y$ then
$$
\hbox{$P'(y) \geq P(y)$ and $P'(x) > P(x)$ $\Rightarrow$}
$$
$$
\hbox{$\Rightarrow$ $\hbar(P'(x),P'(y)) > \hbar(P(x),P(y)) \geq l(G,P)$}.
$$
If $y \in Y$ note that by
$$
\eps \leq \frac{q(t-1)}{(Rt+1)},
$$
one has
$$
1 < \frac{q+\eps}{q-\nu} \leq t.
$$
By Property~\ref{ip2}, setting $\alpha=|Y|$ and $\tau=|\Gamma_FY|$:
$$
\hbar(P(x),P(y))= \hbar(q,q)=z(1,w,\alpha,\tau) < \hbar(P'(x),P'(y)) \leq
z(t,w,\alpha,\tau).
$$

\medskip

{\bf Case 2: }$x \in Y$ and $y \not\in \Gamma_FY$. Clearly, it follows that
$y$ belongs to a component $F' \neq F$. Note that $F' \not\in {\cC}$, otherwise by
definition of ${\cC}$, $q(F')=q$ would imply $p(F')=q$, and $F=F'$. In addition by $x \in
e(P)$ it follows $y \not\in m(P)$. Hence $P(y)=p(F')$, by
$$
\eps \leq \min \{ R^{-1}[q-q(F')] \,,\,[p(F')-q] \}
$$
one has
$$
q(F') \leq q-\nu \leq q + \eps \leq p(F'),
$$
and:
$$
\hbar(P'(x),P'(y))=\hbar(q-\nu ,p(F')) \geq \hbar(q(F'),p(F'))=l(G,P).
$$
Now, note that (Case~1) no nodes in $\Gamma_FY$ are endpoints of edges in ${\cL}(P')$
and
so $P'$ is not $G$-balanced.
\qed

\medskip

Now, we prove item $1$ in Lemma~\ref{main}. For an arbitrary $maximal$
stable
set $X$ such that $P \in Cr(X)$, let us introduce the following relation between the
components of the graph $(V(G),{\cal L}(P))$:
\begin{equation}
\label{prec}
F \prec F'  \:\: \hbox{ iff } \:\:
F \neq F' \: \hbox{ and } \:\exists\, \{x,y\} \in E(G):\, \: x \in X(F), \,
y \in V(F'),
\end{equation}
For the transitive closure ${\stackrel{.}{\prec}}$ of $\prec$ we prove
\begin{claim}
\label{order}
$$
\hbox{$F{\stackrel{.}{\prec}} F'$  $\Rightarrow$ $q(F') < q(F) \leq p(F)<p(F').$}
$$
\end{claim}
\proof
Let $F=F_1 \prec F_2 \prec \ldots \prec F_m=F'$ be any chain of relations
$\prec$. We show
$$
q_m<\ldots<q_2<q_1 \:\leq\: p_1<p_2<\ldots<p_m
$$
where $q(F_i)=q_i$ and $p(F_i)=p_i$. By definition of $\prec$ there exist $m-1$ edges
$(u_i,v_{i+1})$ such that for each $1 \leq i < m$
$$
u_i \in X(F_i) \:\:\hbox{ and }\:\: v_{i+1} \in \bar{X}(F_{i+1})
$$
and
$$
\hbar(P(u_i),P(v_{i+1}))=\hbar(q_i,p_{i+1}) > l(G,P)=\hbar(q_i,p_i),
$$
where the strict inequality follows from $(u_i,v_{i+1}) \not\in {\cL}(P)$. By
Property~\ref{ip1} it follows the claim.
\qed

\medskip

{\flushleft{\bf Observation 2: }}Note that by Claim~\ref{order} if
$q(F)=p(F)$ and $P \in Cr(X)$, then
$$
X(F)=X \cap \Gamma \bar{X}(F),
$$
or else there would exist a component $F' \prec F$.

\medskip

The proof of the following property of $\hbar$ can be found in Appendix~A

\begin{prop}
\label{ip3}
If $\alpha / \tau >1$ then for any $w >0$
$$
\frac{\alpha}{\tau} < \frac{\alpha'}{\tau'} \:\:\hbox{ iff }\:\:
t(w,\tau,\alpha) < t(w,\tau',\alpha').
$$
\end{prop}
Note also that $t$ is independent by $w$. That is, for any $\alpha$, $\tau$,
$w$ and $w'$ (see Appendix~A):
\begin{equation}
\label{w-ind}
t(w,\tau,\alpha)=t(w',\tau,\alpha)=t(\tau,\alpha)
\end{equation}
\begin{lem}
\label{GEQ}
If $P$ is a $G$-balanced distribution centered on a stable set $X$ and
$F$ is a connected component of the graph $(V(G),{\cal L}(P))$, then for
any $U \subseteq \bar{X}(F)$
$$
t(|X(F)\cap \Gamma U |,|U|) \geq \frac{p(F)}{q(F)}.
$$
\end{lem}
\proof
Suppose that the above inequality is false for $U \subseteq \bar{X}(F)$. As in the proof
of
item $2$, we replace $P$ with a new probability distribution $P'$ where
${\cL}(P')$ is not a line cover of $G$ and $l(G,P') \geq l(G,P)$.

Set
$$
t=t(|X(F)\cap \Gamma U |,|U|).
$$
By hypothesis
$$
t<\frac{p(F)}{q(F)}.
$$
Fix
$$
\hbox{$R=\frac{|U|}{|X(F)\cap \Gamma U|}$   and   $L=\{F':\,F'{\stackrel{.}{\prec}}F \}$}
$$
and the two real numbers
$$
\epsilon = \min\left\{ \frac{p(F)-q(F)t}{Rt+1}
\,,\,R^{-1}\min_{F' \in L}{[q(F')-q(F)]} \,,\,\min_{F' \in
L}{[p(F)-p(F')]} \right\},
$$
$$
\nu = \epsilon R,
$$
by Claim~\ref{order} and $p(F)>q(F)t$, $\eps>0$. Define $P'$ as:
\begin{equation}
P'(v)=\cases{p(F)-\epsilon & if $v\in U$ \cr
                    & \cr
            q(F)+\nu & if $v\in X(F) \cap \Gamma U$  \cr}
\end{equation}
\medskip
and $P'(v)=P(v)$ if $v \not\in U \cup (X(F) \cap \Gamma U)$.

$P'$ is a probability distribution. Indeed:
$$
p(F)-\epsilon \geq \frac{t(q(F)+Rp(F))}{Rt+1}>0,
$$
and $\nu$ is chosen so as to leave unchanged the total amount of probability of $U \cup
(X(F) \cap \Gamma U)$. Also note that by $t \geq 1$,
$$
p(F)-\eps \geq q(F)+\nu,
$$
that is
$$
\eps \leq \frac{p(F)-q(F)t}{Rt+1} \leq \frac{p(F)-q(F)}{R+1}.
$$
Now, we show that $l(G,P') \geq l(G,P)$.

\medskip

{\bf Case 1: }Edges $\{x,y\}$ such that $x \in X(F) \cap \Gamma U$.
If $y \not\in U$ then $P'(y) \geq P(y)$ and $P'(x) > P(x)$ make this case trivial. If $y
\in
U$, by
$$
\eps \leq \frac{(p(F)-q(F)t)}{(Rt+1)},
$$
one has
$$
\frac{p(F)}{q(F)} > \frac{p(F)-\eps}{q(F)+\nu} \geq t.
$$
Now, set
$$
\hbox{$\alpha=|X(F) \cap \Gamma U|$, $\tau=|U|$}
$$
and
$$
P(F)=\sum_{v \in V(F)}P(v),
$$
By Property~\ref{ip2}, and~(\ref{w-ind}) one has
$$
\hbar(P(x),P(y))= \hbar(p(F),q(F))=z \left( \frac{p(F)}{q(F)},P(F),\alpha,\tau \right)
$$
$$
<\hbar(P'(x),P'(y)) \leq  z(t,P(F),\alpha,\tau).
$$

\medskip

{\bf Case 2: }$x \in U$ and $y \not\in X(F) \cap \Gamma U$. Let $y \in V(F')$. If $y \in
\bar{X}(F')$ suppose $p(F') \geq p(F)-\eps$, then
$$
\hbar(P'(x),P'(y))=\hbar(p(F)-\eps,p(F')) \geq \hbar(p(F)-\eps,p(F)-\eps)\geq
$$
$$
\geq \hbar(p(F)-\eps,q(F)+\nu)>\hbar(p(F),q(F))=l(G,P).
$$
The case $p(F') < p(F)-\eps$ can be evaluated in a similar way.

Finally, if $y \in X(F')$ then by hypothesis $F' \neq F$ and by
$U \subseteq \bar{X}(F)$, $F' \prec F$. Being
$$
\eps \leq \min \left \{ R^{-1}[q(F')-q(F)] \,,\,[p(F)-p(F')] \right \},
$$
it follows
$$
q(F)+\nu \leq q(F') \leq p(F') \leq p(F)- \eps.
$$
That is
$$
\hbar(P'(x),P'(y))=\hbar(p(F)-\eps,q(F')) \geq
$$
$$
\geq \hbar(p(F)-\eps,q(F)+\nu)>\hbar(p(F),q(F))=
l(G,P).
$$
Now, note that the value of $\hbar(.,.)$ strictly increases over all the edges with at last
one
point in $X(F) \cap \Gamma U$ (Case~1). Hence, unless $U=\bar{X}$ (in this
case one should have directly $t=p(F)/q(F)$) it follows that ${\cL}(P')$ is not a line
cover
of $G$ which proves the statement.
\qed

{\flushleft {\bf Proof of $1$ in Lemma~\ref{main}:}}
Let ${\cal I}(.)$ be the family of all the $maximum$ stable sets in a graph.
Note that $1$ is equivalent to the following equality
\begin{equation}
\label{stable-eq}
{\cal I}(G)=\{Z:\,Z=m(P)\cup A,\,A\in {\cal I}(F)\}.
\end{equation}

Let us consider any maximum stable set $Z$ in $G$. Then, if we show that for every $G$-
balanced distribution $P$, $P \in Cr(Z)$ this would imply~(\ref{stable-eq}). Indeed, by
definition of $Cr(Z)$,
$$
Z=A \cup m(P),
$$
where $A$ is a maximal stable set in the subgraph induced in $G$ by $e(P)$. By
Claim~\ref{disjoint}, if $S$ is any stable set in such a subgraph then the set $Z'=S \cup
m(P)$ is a stable set in $G$. Hence,
$$
\hbox{$|Z'|=|S|+|m(P)| \leq |Z|=|A|+|m(P)|$ $\Rightarrow$ $A \in {\cal I}(F)$}.
$$
Vice versa, once again by Claim~\ref{disjoint}, if $S \in {\cal I}(F)$ then $S \cup m(P)$
is a
maximal stable set in $G$ and by $S$ maximum in $F$
$$
|S|+|m(P)| \geq |A|+|m(P)|.
$$
On the other hand, we supposed $Z$ maximum that is $S \cup m(P) \in {\cal I}(G)$.

Suppose $P \in Cr(Y)$. If $Z=Y$ we have finished. Let $F$ be any connected component
of
the graph $(V(G),{\cL}(P))$, and set
$$
\Delta_{Z}(F)=V(F) \cap (Z \setminus Y).
$$
Let us fix
$$
{\cC}=\{F:\, \Delta_{Z}(F) \neq \emptyset \}.
$$
Being ${\cL}(P)$ a line cover, $\{\Delta_{Z}(F):\,F \in \cC \}$ is a partition of $Z
\setminus
Y$.  Further, if $F \in {\cC}$ then
$$
\Delta_Y(F)=V(F) \cap (Y \setminus Z) \neq \emptyset.
$$
Indeed, if $x \in \Delta_{Z}(F)$ then $x \in \bar{Y}(F)$ and, being ${\cal L}(P)$ a line
cover of $G$, $Y(F)\cap \Gamma x  \neq \emptyset$. Now, by $x \in Z \setminus Y$ it
follows
$$
Y(F) \cap \Gamma x \subseteq Y \setminus Z.
$$
So, being
$$
|Z \setminus Y|= \left | \bigcup_{F \in {\cC}}\Delta_{Z}(F) \right | = \sum_{F \in
{\cC}}|\Delta_{Z}(F)|
$$
and
$$
|Y \setminus Z| \geq \left | \bigcup_{F \in {\cC}}\Delta_Y(F) \right |
=\sum_{F \in {\cC}}|\Delta_Y(F)|,
$$
we have
$$
\min_{F \in {\cC}}\frac{|\Delta_Y(F)|}{|\Delta_Z(F)|} \leq
\frac{\sum_{F \in {\cC}}|\Delta_Y(F)|}{ \sum_{F \in {\cC}}|\Delta_Z(F)|} \leq \frac{|Y
\setminus Z|}{|Z \setminus Y|} \leq 1.
$$

\medskip

Therefore, we can fix any $C \in {\cC}$ such that
$$
|\Delta_Y(C)| \leq |\Delta_Z(C)|.
$$
By $\Delta_Z(C) \subseteq Z \setminus Y$ it follows
$Y~\cap~\Gamma~[\Delta_Z(C)]~\subseteq~Y~\setminus~Z$ and
in particular
$$
Y(C) \cap \Gamma [\Delta_Z(C)] \subseteq \Delta_Y(C).
$$
Hence
\begin{equation}
\label{re1}
\frac{|Y(C) \cap \Gamma [\Delta_Z(C)]|}{|\Delta_Z(C)|} \leq
\frac{|\Delta_Y(C)|}{|\Delta_Z(C)|} \leq 1.
\end{equation}
In accordance with Lemma~\ref{GEQ} and~(\ref{down})
$$
\frac{p(C)}{q(C)} \leq t(|Y(C) \cap \Gamma [\Delta_Z(C)]|,|\Delta_Z(C)|)=1.
$$
So $p(C)=q(C)$ and $\Delta_Z(C) \subseteq e(P)$. Moreover, by Observation~2 it follows
that
$Y(C)=Y \cap \Gamma \bar{Y}(C)$ and then
\begin{equation}
\label{loc}
Y \cap \Gamma [\Delta_Z(C)]=Y(C) \cap \Gamma [\Delta_Z(C)] \subseteq e(P).
\end{equation}

Now, set
$$
K=(Y \setminus \Gamma [\Delta_Z(C)]) \cup \Delta_Z(C).
$$
Note that by~(\ref{re1}) and~(\ref{loc}), $|K| \geq |Y|$ and it is easy to check that $K$ is
a stable set in $G$. By
$$
\Delta_Z(C)\cup(Y \cap \Gamma [\Delta_Z(C)]) \subseteq e(P),
$$
and $m(P) \subseteq Y$, it follows $m(P) \subseteq K$.
If $R_K$ is any maximal stable set in the subgraph induced by $e(P)$ in $G$ containing
$K
\setminus m(P)$ then by Claim~\ref{disjoint}
$$
Y_1=R_K \cup m(P)
$$
is a stable set in $G$. We have $|Y_1| \geq |K| \geq |Y|$ and $P \in Cr(Y_1)$. In addition
$$
|Z \setminus Y_1| \leq |Z \setminus K|=|Z \setminus Y|-|\Delta_Z(C)|<|Z \setminus Y|,
$$
and
$$
\hbox{$|Z \setminus Y_1|=0$ implies $Y_1=Z$}.
$$
Iteratively applying the above procedure, we find a sequence of maximal stable
sets: $Y=Y_0,Y_1,\ldots$ such that $P \in Cr(Y_i)$ and $|Z \setminus Y_i|$ $strictly$
decreases with $i \geq 0$. Hence, for some $m>0$ we get $Y_m=Z$ and the statement.
\qed

\section{Appendix A: basic properties of $\hbar$}

We prove the three main properties of $\hbar$. Property~\ref{ip1} is easy to verify.
\begin{enumerate}
\item{\bf Property~\ref{ip2}, (\ref{w-ind}), (\ref{down})} We have:
\begin{equation}
z(t,w,\alpha,\tau)=\frac{w}{t\tau+\alpha}\left[ \log(t+1)+\log \left(\,
\frac{1}{t}+1 \right) \right]
\end{equation}
and:
$$
\frac{dz}{dt}=\frac{w}{(t\tau+\alpha)^2}\left[\alpha\log \left(\,
\frac{1}{t}+1 \right)-\tau\log(t+1)\right].
$$
Hence $t(w,\alpha,\tau)$ is independent by $w$ and it follows~(\ref{w-ind}).
Now, if $\alpha \leq \tau$ the point of maximum of $z$ is $t(\tau,\alpha)=1$. Otherwise
$t(\tau,\alpha)$ is the unique number greater than 1 that is a
root of:
$$
\rho(t)=(t+1)^{\alpha-\tau}-t^\alpha.
$$
This proves Property~\ref{ip2} and~(\ref{down}). We note that $(\ref{down})$
can be proved for any function verifying Properties~\ref{ip1} and~\ref{ip2} (the
proof is
not trivial).

\item{\bf Property~\ref{ip3}}: Remember that if $\alpha > \tau$, $t=t(\alpha,\tau)$ is
the unique root greater than $1$ of
$$
\rho(t)=(t+1)^{\alpha-\tau}-t^\alpha.
$$
Hence
$$
\frac{\tau}{\alpha}=1-\frac{\log t}{\log (t+1)}
$$
and it is sufficient to note that the right hand side is a strictly decreasing function on
the semi-interval $t \geq 1$.
\end{enumerate}

\message{References}


\begin{thebibliography}{10}

\bibitem{BS}
M.L. Balinski and K.~Spielberg.
\newblock Methods for integer programming: algebraic, combinatorial and
  enumerative.
\newblock In J.~Aronofsky, editor, {\em Progress in operation research}, volume
  III, pages 195--292, Wiley, New York, 1969.

\bibitem{Ber}
C.~Berge.
\newblock {\em Graphes}.
\newblock Gaulhier-Villars, Paris, 1973.

\bibitem{CKS}
G.~Cohen, J.~K{\"o}rner, and G.~Simonyi.
\newblock Zero error capacities and very different sequences.
\newblock In R.M. Capocelli, editor, {\em Sequences: combinatorics, compression
  security and trasmission}, pages 144--155. Springer-Verlag, 1990.

\bibitem{GKV1}
L.~Gargano, J.~K{\"o}rner, and U.~Vaccaro.
\newblock Sperner capacities.
\newblock {\em Graphs and combinatorics}, 9:31--46, 1993.

\bibitem{Gre95}
G.~Greco.
\newblock Capacities of graphs and $2$-matchings.
\newblock {\em Discrete Mathematics}, 186:135--143, 1998.

\bibitem{Lor}
L.C. Lorentzen.
\newblock {\em Notes on covering of arcs by nodes in an undirected graph}.
\newblock Technical report, 1966.

\bibitem{LP}
L.~Lov{\'a}sz and M.D. Plummer.
\newblock {\em Matching Theory}.
\newblock North-Holland, New-York, 1986.

\bibitem{Pul}
W.R. Pulleyblank.
\newblock Minimum node covers and $2$-bicritical graphs.
\newblock {\em Mathematical programming}, 17:91--103, 1979.

\bibitem{Re}
A.~R{\'e}nyi.
\newblock {\em Probability theory}.
\newblock North-Holland, Amsterdam/New-York, 1970.

\bibitem{Tut}
W.T. Tutte.
\newblock The $1$-factors in oriented graphs.
\newblock {\em Proceedings american mathematical society}, 22:107--111, 1947.

\end{thebibliography}
\end{document}